\documentclass[final,5p,authoryear,times]{elsarticle}
\usepackage[T1]{fontenc}
\pdfoutput=1
\usepackage[colorlinks,linkcolor=blue,citecolor=blue,urlcolor=blue]{hyperref}
\usepackage{graphicx}
\usepackage{amsfonts,amsmath,amssymb}
\usepackage{mathptmx}
\usepackage{bm}
\usepackage{dblfloatfix} 
\usepackage{caption}
\usepackage{booktabs}
\usepackage[ruled]{algorithm}
\usepackage[noend]{algpseudocode}
\alglanguage{pseudocode}
\algrenewcommand{\algorithmiccomment}[1]{\hskip3em\# #1}

\usepackage{xspace}

\newcommand{\tens}[1]{\mathsf{#1}}
\renewcommand*{\vec}[1]{\mathbf{#1}}
\newcommand{\x}{\vec{x}}
\newcommand{\z}{\vec{z}}
\newcommand{\m}{\vec{m}}
\renewcommand{\v}{\vec{v}}
\newcommand{\g}{\vec{g}}
\renewcommand{\u}{\vec{u}} 
\newcommand{\tA}{\tens{A}}
\newcommand{\tS}{\tens{S}}

\newcommand{\tY}{\tens{Y}}
\newcommand{\tP}{\tens{P}}
\newcommand{\tH}{\tens{H}}
\newcommand{\tD}{\tens{D}}
\newcommand{\prox}{\mathrm{prox}}
\providecommand{\norm}[1]{\lVert#1\rVert}
\newcommand{\adagrad}{{\sc Ada\-Grad}\xspace}
\newcommand{\rmsprop}{{\sc RMS\-Prop}\xspace}
\newcommand{\amsgrad}{{\sc AMS\-Grad}\xspace}
\newcommand{\adam}{{\sc Adam}\xspace}
\newcommand{\adamx}{{\sc AdamX}\xspace}
\newcommand{\padam}{{\sc PAdam}\xspace}
\newcommand{\adaprox}{{\sc Ada\-Prox}\xspace}

\newcommand{\software}[1]{{\tt #1}\xspace}

\begin{document}
\begin{frontmatter}

\title{Proximal Adam: Robust Adaptive Update Scheme for Constrained Optimization}
\author[ad1,ad2]{Peter Melchior}
\ead{peter.melchior@princeton.edu}
\author[ad1]{R\'emy Joseph}
\author[ad1]{Fred Moolekamp}
\address[ad1]{Department of Astrophysical Sciences, Princeton University, Princeton, NJ 08544, USA}
\address[ad2]{Center for Statistics \& Machine Learning, Princeton University, Princeton, NJ 08544, USA}

\begin{abstract}
We implement the adaptive step size scheme from the optimization methods \adagrad and \adam in a novel variant of the Proximal Gradient Method (PGM).
Our algorithm, dubbed \adaprox, avoids the need for explicit computation of the Lipschitz constants or additional line searches and thus reduces per-iteration cost. 
In test cases for Constrained Matrix Factorization we demonstrate the advantages of \adaprox in fidelity and performance over PGM, while still allowing for arbitrary penalty functions.
The {\sf python} implementation of the algorithm presented here is available as an open-source package at \url{https://github.com/pmelchior/proxmin}.
\end{abstract}


\end{frontmatter}

\section{Introduction}

Many problems in design, control, and parameter estimation seek to
\begin{equation}
\label{eq:problem}
\underset{x\in\mathbb{R}^d}{\mathrm{minimize}}\  f(\x) + r(\x),
\end{equation}
where $f$ is a smooth convex loss function with a Lipschitz-continuous gradient, i.e. 
\begin{equation}
\exists\, L: \norm{\nabla f(\x) - \nabla f(\z)}\leq L\norm{\x-\z},
\end{equation}
and $r$ is a convex, potentially non-differentiable penalty function that regularizes the solution.
For instance, for parameter estimation $f$ is the log-likelihood of $\x$ given some observations, and $r$ represents a predetermined solution manifold or subspace.

First-order gradient methods usually depend in some form on $L$ to determine the size of gradient steps $\alpha \propto 1/L$.
We seek to find an algorithm that avoids the computation of $L$, which can be costly in practice, while permitting arbitrary regularizers $r$, as long as they can be expressed through their proximal operators \citep{Moreau-1965}
\begin{equation}
\label{eq:proximal}
\prox_{\alpha\, r}(\x)\equiv\underset{\z}{\textrm{argmin}}\left\lbrace r(\z)+\frac{1}{2\alpha}\norm{\z-\x}_{2}^{2}\right\rbrace.
\end{equation}
Motivated by data-intensive analysis problems, we are particularly concerned with problems for which $f$ and its derivatives are expensive to evaluate, but $\prox_r$ is not.
To avoid explicit or implicit computation of $L$, we instead employ a class of optimization algorithms popularized by deep learning applications: \adam and its variants \amsgrad, \adamx, \padam.

This paper is structured as follows. \autoref{sec:opt} provides the motivation for this work and reviews the relevant proximal and adaptive optimization techniques. \autoref{sec:adaprox} introduces our new adaptive proximal method \adaprox. \autoref{sec:cmf} compares PGM with \adaprox on three variants of constrained matrix factorization. We conclude in \autoref{sec:conclusions}.

\section{Proximal and Adaptive Optimization}
\label{sec:opt}

\begin{table*}[h]
\caption{Choices to accumulate mean and variance of $\g\equiv\nabla f(\x)$ for the algorithms discussed in this work, typically via intermediate variables $\m_t$, $\v_t$, and $\hat\v_t$.
Steps sizes $\alpha_t$ are usually set to $\alpha/\sqrt{t}$ for provable convergence, but in practice often follow a different schedule, including constant steps. PGM uses $\alpha_t \in(0, 2/L_t)$, usually $1/L_t$. Constants $\beta_1$ or scheduled $\beta_{1,t}$, $\beta_2$ are from $[0,1)$, $\epsilon >0$, and $p \in (0,1/2]$.}
\label{tab1}
\begin{tabular}{lcccccl}
\toprule
& \multicolumn{2}{c}{Mean estimate} & \multicolumn{3}{c}{Variance estimate}\\
 \cmidrule(r){2-3} \cmidrule(r){4-6}
Name & $\m_t$ & $\phi_t$ & $\v_t$ & $\hat{\v}_t$ & $\psi_t$\\
\midrule
SGD, PGM & --- & $\g_t$ & --- & --- & $\mathbb{I}$ \\

\adagrad & --- & $\g_t$ & ---  & --- & $\sqrt{\frac{1}{t}\sum_{i=1}^t \g_i^2}$ \\

\adam & $\beta_1 \m_{t-1} + (1-\beta_1)\g_t$ & $\m_t / (1 - \beta_1^t)$  & $\beta_2 \v_{t-1} + (1-\beta_2)\g_t^2$ & --- & $\sqrt{\v_t / (1 - \beta_2^t)} + \epsilon$ \\

\amsgrad & $\beta_{1,t} \m_{t-1} + (1-\beta_{1,t})\g_t$ & $\m_t$ &  $\beta_2 \v_{t-1} + (1-\beta_2)\g_t^2$ & $\max\left(\hat\v_{t-1}, \v_t\right)$ & $\sqrt{\hat\v_t}$ \\

\adamx & $\beta_{1,t} \m_{t-1} + (1-\beta_{1,t})\g_t$ & $\m_t$ &  $\beta_2 \v_{t-1} + (1-\beta_2)\g_t^2$ & $\max\left(\frac{(1-\beta_{1,t})^2}{(1-\beta_{1,t-1})^2}\hat\v_{t-1}, \v_t\right)$ & $\sqrt{\hat\v_t}$ \\

\padam & $\beta_{1,t} \m_{t-1} + (1-\beta_{1,t})\g_t$ & $\m_t$ &  $\beta_2 \v_{t-1} + (1-\beta_2)\g_t^2$ & $\max\left(\hat\v_{t-1}, \v_t\right)$ & $\hat\v_t^p$ \\
\bottomrule
\end{tabular}
\end{table*}

\subsection{Proximal Gradient Method}
A well-known and effective approach for solving \autoref{eq:problem} is a \emph{forward-backward} scheme, where at iteration $t$ a step in the direction of $\nabla f$ is followed by the application of the proximal operator:
\begin{equation}
\label{eq:pgm}
\x_{t+1} = \prox_{\alpha_t r} \left(\x_t - \alpha_t \nabla f(\x_t)\right).
\end{equation}
If step size $\alpha_t \in (0, 2/L)$, the sequence converges to the minimum of $f+r$.
This algorithm is known as Proximal Gradient Method \citep[PGM, e.g.][]{Parikh2014-ex}.

It is straightforward to compute the Lipschitz constants for simple problems, but more complex problems can make that computation non-analytic or very expensive.
As an example, consider the linear inverse problem,
\begin{equation}
f(\x) = \frac{1}{2}\norm{\tP\x - \vec{y}}_2^2
\end{equation}
for some observation $\vec{y}$ with i.i.d. Gaussian errors.
The gradients of $f$ are bound by $L = \norm{\tP^\top \tP}_s$\footnote{In this work, we denote the element-wise 2-norm and the spectral norm, i.e. the largest eigenvalue, as $\norm{.}_2$ and $\norm{.}_s$, respectively.}.
In image analysis the matrix $\tP$ typically encodes resampling and convolution operations, so that $L$ is expensive to compute.
Once the problem includes non-linear mappings, e.g.
\begin{equation}
\label{eq:non-linear-f}
f(\x) = \frac{1}{2}\norm{\tP s(\x) - \vec{y}}_2^2,
\end{equation}
with some differentiable parameterization $s$ of the signal, $L$ becomes a function of $\x$ and has to be recomputed at each iteration of the optimizer.
That is also true in multi-convex cases such as matrix or tensor factorization.
To make matters worse, $L$ often does not have an analytically known form and thus has to be determined by yet another procedure like a line search \citep{Beck2009-bf}.
While effective, it requires multiple evaluations of $f$ per optimization parameter, which quickly becomes prohibitive if $f$ is expensive to evaluate.

We therefore seek a proximal gradient method whose step sizes can be set without invoking Lipschitz constants, but maintain the flexibility of PGM to accept arbitrary regularizers. 
For instance in astronomy, the regularization is almost always physically motivated and can drastically vary between different analyses.
One could also avoid the limitations of PGM with second-order methods in the form of a proximal (quasi-)Newton scheme \citep[e.g.][]{Becker2012-ry, Becker2019-mv}. It replace stepsizes $\propto 1/L$ by multiplications with the inverse Hessian matrix of $f$, but computing the Hessian is at least as expensive as computing $L$.

\subsection{Adaptive gradient methods}

In machine-learning applications such as the training of deep neural networks, generic optimizers are routinely employed for any functional form of the model or the loss function.
This makes it hard to determine reasonable step sizes \emph{a priori}, and the problem sizes are often too large to allow for line searches on-the-fly.
An adaptive optimizer that excels in such cases is \adam\ \citep{Kingma2015-pq}.

The central idea for adaptive gradient updates amounts to replacing a simple gradient step with
\begin{equation}
\label{eq:adaptive}
\x_{t+1}=\x_t - \alpha_t \frac{\m_t}{\sqrt{\v_t}}
\end{equation}
where $\alpha_t$ are externally provided step sizes, potentially varying at every step $t$; and $\m_t$ and $\v_t$ are estimates of the mean and variance of $\g\equiv\nabla f(\x)$, respectively.
This scheme has two effects: 1) It adjusts the step size for every dimension as a cheap emulation of a Newton scheme. 2) It renders the updates steps $\alpha_t$ independent of the actual amplitude of $\g$, sidestepping the problem of having to compute a Lipschitz constant.
Moreover, for physically motivated problems it is more natural to think of step sizes in the units of the parameter instead of the units of $f$.

\adagrad\ \citep{Duchi2011-bb}, one of the first algorithms to use the scheme, was designed for online optimization with sparse gradients and therefore sums up $\g^2$ from all previous iterations as $\v_t$.
\rmsprop\ \citep{Hinton2012-rx} and \adam\ maintain the general form of \autoref{eq:adaptive} but replace the moment accumulation with exponential moving averages, which has proven very successful in practice, especially for stochastic gradients. 
More recently, flaws in the original convergence proof of \citet{Kingma2015-pq} have triggered a series of minor modifications to the form of the $\v_t$ term, e.g. \amsgrad\ \citep{Reddi2018-me}, \padam\ \citep{Chen2018-wu}, and \adamx\ \citep{Phuong2019-eh}.
To better clarify the corresponding choices, we rewrite \autoref{eq:adaptive} as
\begin{equation}
\label{eq:ada-phi-psi}
\x_{t+1}=\x_t - \alpha_t \frac{\phi(\g_1,\dots,\g_t)}{\psi(\g_1^2,\dots,\g_t^2)}
\end{equation}
and list the choices for $\phi$ and $\psi$ of each algorithm in \autoref{tab1}.

\section{Adaptive proximal gradient methods}
\label{sec:adaprox}

We seek to combine the general purpose forward-backward splitting method of \autoref{eq:pgm} with the robustness and efficiency of the adaptive gradient update scheme of \autoref{eq:ada-phi-psi}.
The introduction of $\psi$ updates every dimension $j$ of $\x$ with a different effective learning rate $\alpha_t / \psi_{t,j}$, which is equivalent to introducing a metric $\tH_t$ for the parameter space.
If $\psi$ is an approximation of the Hessian of $f$, the update corresponds to a proximal quasi-Newton scheme \citep{Becker2012-ry,Tran-Dinh2015-bp} of the form
\begin{equation}
\label{eq:ada-pgm}
\x_{t+1} = \prox_{\alpha_t r}^{\tH_t} \left(\x_t - \alpha_t \frac{\phi(\g_1,\dots,\g_t)}{\psi(\g_1^2,\dots,\g_t^2)} \right),
\end{equation}
with a variable-metric proximal operator
\begin{equation}
\label{eq:scaled-prox}
\prox_{\alpha r}^\tH(\x)\equiv \underset{\z}{\textrm{argmin}}\left\lbrace r(\z)+\frac{1}{2\alpha}\norm{\z-\x}_{\tH}^{2}\right\rbrace
\end{equation}
and  $\norm{\x}_\tH^2 \equiv \x^\top \tH \x$.
Were one to apply the regular proximal operator in an adaptive scheme without considering the variable metric, the results would be feasible but not optimal.

The metric $\tH_t$ does not need to approximate the Hessian of $f$.
\adagrad \citep{Duchi2011-bb} introduced a variable-metric projection $\Pi^\tH_\mathcal{S}$ of the updated parameter $\x_{t+1}$ onto a convex subset $\mathcal{S}\subset\mathbb{R}^d$ with $\phi_t / \psi_t$ from \autoref{tab1}.
In particular, $\phi_t=\g_t$, i.e. the instantaneous gradient direction.
Later methods have adopted moving averages, equivalent to an inexact proximal gradient method, which does not affect convergence as long as $\phi_t - \g_t$ decreases as $\mathcal{O}\left(1/t^{1+\delta}\right)$ for any $\delta>0$ \citep{Schmidt2011-lu}.
\adagrad is limited to projection operators, which are a special class of proximal operators, namely those for the indicator function of any convex subset $\mathcal{S}\subset\mathbb{R}^d$.
A limitation to indicator functions is not fundamental, by lifting it we seek to allow for regularizers that impose e.g. sparsity or low-rankness of the solutions.

This brings us back to the question how to solve \autoref{eq:scaled-prox}.
\citet{Chouzenoux2014-wg} proposed a dual forward-backward algorithm from \citet{Combettes2011-dd}, which requires computing $L$ and is thus not applicable here.
\citet{Becker2012-ry} showed that if $\tH \equiv \tD + \u \u^\top$ with a diagonal $\tD$ and an arbitrary $\u\in\mathbb{R}^d$, $\prox_{\alpha r}^\tH(\x)$ can be replaced with $\prox_{\alpha\, r\circ \tD^{-1/2}}(\tD^{1/2}\x - \v)$.
The offset $\v$ needs to be found through a line search that involves $\prox_{\alpha\, r\circ \tD^{-1/2}}$, which itself may be expensive to compute even if $\prox_{\alpha\, r}$ is efficient.
\citet{Becker2019-mv} demonstrated how to perform this computation more directly for several classes of common regularizers, which can lead to substantial performance gains.

We propose a more direct approach that is entirely agnostic about the regularizer.
Because the $\tH$-norm part of \autoref{eq:scaled-prox} is differentiable with gradient $\tfrac{1}{\alpha} \tH (\z- \x)$ and Lipschitz constant $L_\tH=\tfrac{1}{\alpha} \sqrt{\norm{\tH^\top \tH}_s}$, the minimizer of \autoref{eq:scaled-prox} for a given $\x$ can be found with PGM:
\begin{equation}
\z_{\tau+1} = \prox_{\gamma\, r}\left(\z_\tau - \frac{\gamma}{\alpha} \tH(\z_\tau -\x)\right) \ \mathrm{for}\ \tau = 1,2,\dots
\end{equation}
The step size of the sub-problem is as usual $\gamma\in(0, 2/ L_\tH)$.
\citet{Duchi2011-bb} and \citet{Tran-Dinh2015-bp} showed that it is often sufficient, and much more efficient in high-dimensional settings, to diagonalize the metric: $\tH_t = \mathrm{Diag}(\psi_t)$.
With corresponding step sizes $\gamma_t = \alpha_t / \max (\psi_t)$, the proximal sub-iteration to achieve optimality is
\begin{equation}
\label{eq:subiter}
\z_{\tau+1} = \prox_{\gamma_t\, r}\left(\z_\tau - \frac{1}{\max(\psi_t)} \mathrm{Diag}(\psi_t) (\z_\tau-\hat\x_{t+1})\right),
\end{equation}
where $\hat\x_{t+1}$ denotes the unconstrained parameter after gradient update from \autoref{eq:ada-phi-psi}.
Once the desired level of convergence of the $\z$-sequence is reached, $\x_{t+1} \leftarrow \z_{\tau+1}$.
In essence, the PGM sub-iterations enable ordinary proximal operators to be used instead of variable-metric operators that arise from adaptive updates.
The entire algorithm is listed as \autoref{alg1}.

\begin{algorithm}[t]
\caption{Adaptive Proximal Gradient Method (\adaprox)\newline
\small  The constrained problem of minimizing $f+r$ is solved by gradient decent with an adaptive scheme from \autoref{tab1} followed by the solution of the scaled proximal operator (lines 10-12).}
\label{alg1}
\begin{algorithmic}[1]
\Procedure{AdaProx}{$\x_1;\nabla f(.); \prox_r(.); \lbrace\alpha_t\rbrace_t; \lbrace\beta_{1,t}\rbrace_t; \beta_2; \epsilon$}
\For{$t = 1,2,\dots$}
	\State $\g_t = \nabla f(\x_t)$
	\State $\phi_t = \phi(\g_1,\dots,\g_t; \beta_{1,t})$
	\State $\psi_t = \psi(\g_1^2,\dots,\g_t^2; \beta_2)$
	\State $\hat\x_{t+1} = \x_t - \alpha_t \phi_t / \psi_t$
	\State $\tH_t = \mathrm{Diag}(\psi_t)$
	\State $\gamma_t = \alpha_t / \max (\psi_t)$
	\State $\z_{1} = \hat\x_{t+1}$
	\For{$\tau = 1,2,\dots$}
		\State $\z_{\tau+1} = \prox_{\gamma_t r}\left(\z_\tau - \frac{\gamma_t}{\alpha_t} \tH_t (\z_\tau-\hat\x_{t+1})\right)$
		\If{$\norm{\z_{\tau+1} - \z_\tau} < \epsilon \norm{\z_{\tau+1}}$}
			break
		\EndIf
	\EndFor
	\State{$\x_{t+1} = \z_{\tau+1}$}
	\If{$\norm{\x_{t+1} - \x_t} < \epsilon \norm{\x_{t+1}}$}
		break
	\EndIf
\EndFor
\EndProcedure
\end{algorithmic}
\end{algorithm}

\section{Applications to Constrained Matrix Factorization}
\label{sec:cmf}

Matrix Factorization is a non-parametric method for representing a high-dimensional data set through lower-dimensional factors.
In astronomy, applications range from spectral classification to hyperspectral unmixing and multi-band source separation.
We are particularly interested in source separation problems, i.e. we seek to find factors $\tA$ and $\tS$ to approximate an observation matrix $\tY$ by minimizing the loss
\begin{equation}
\label{eq:mf}
f(\x) = \frac{1}{2}\norm{\tA \tS - \tY}_2^2
\end{equation}
with respect to the parameters $\tA$ and $\tS$.
The inherent degeneracies demand additional constraints to be placed on the factors, which requires the use of constrained optimization techniques.
PGM can be used for this problem in an alternating approach of updating $\tA$ at fixed $\tS$ and then $\tS$ at fixed $\tA$ \citep{Rapin2013-xu, Xu2013-nh}.

For this bilinear problem the Lipschitz constants $\norm{\tA_t^\top \tA_t}_s$ and  $\norm{\tS_t \tS_t^\top}_s$ have to be recomputed at every iteration $t$.
It becomes more complicated if data are affected by heteroscedastic or correlated noise. The loss function generalizes to
\begin{equation}
f(\x) = \frac{1}{2}\left(\tA\tS - \tY\right)^\top \Sigma^{-1}  \left(\tA\tS - \tY\right)
\end{equation}
with an inverse covariance matrix $\Sigma^{-1}$.
As we have shown \citep[section 2]{Melchior2018-hb}, the Lipschitz constants remain analytic but involve products of block-diagonal representations of $\tA$ and $\tS$ with $\Sigma^{-1}$, which require spectral norms for very large matrices.
A similar complication arises in online optimization or data fusion applications because not every batch or data set $\tY_l$ has an equal amount of information on all parameters.
The joint loss function for matrix factorization
\begin{equation}
\label{eq:fusion-mf}
f(\x) = \frac{1}{2}\sum_l\left(\tP_l\tA\tS - \tY_l\right)^\top \Sigma_l^{-1}  \left(\tP_l\tA\tS - \tY_l\right)
\end{equation}
has gradients with complicated structure depending on the degradation operators $\tP_l$ and noise properties $\Sigma_l$ of observation $l$.
In particular, the na\"ive estimate $L = \sum_l L_l$ is an upper bound for the joint Lipschitz constant that is applicable only in the unrealistic case that the data sets provide identical information about the parameters.
The resulting step sizes will be under-estimated and thus slow down the convergence of the optimization.
These applications should therefore benefit from our proposed adaptive proximal scheme.

\subsection{Non-negative and Mixture Matrix Factorization}
\label{sec:nmf}

\begin{figure}
\includegraphics[width=\linewidth]{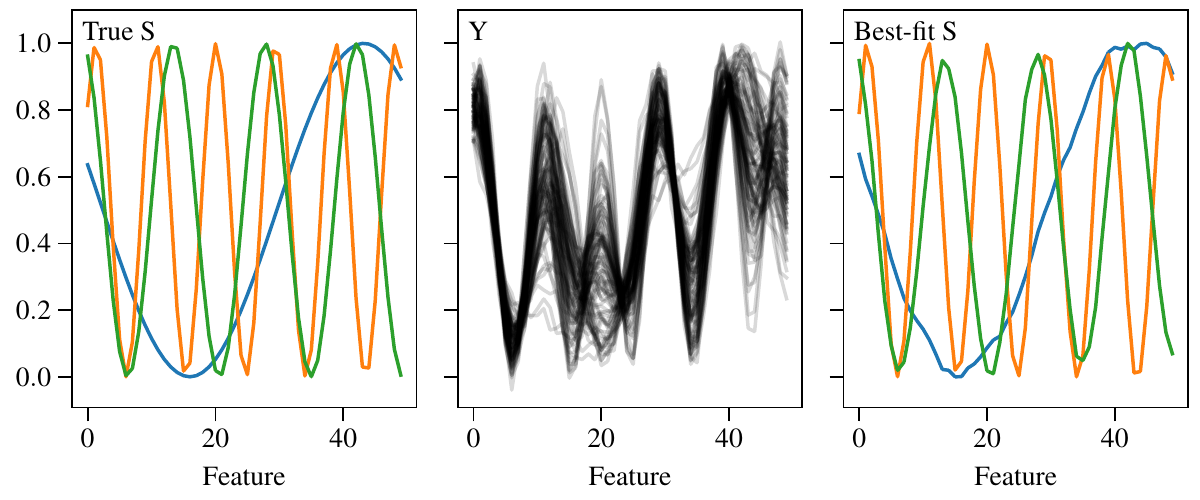}
\caption{NMF test data of $K=3$ sinosoidal components $\tS_k\in\mathbb{R}^{50}$ (\emph{left}), observed 100 times with different mixing weights and i.i.d. Gaussian noise of $\sigma_n = 0.02$ (\emph{center}). The best-fit result of \adaprox-\amsgrad\ with $\alpha=0.1$, rescaled to a maximum of 1, is shown in the \emph{right} panel. The test data are publicly available in the code repository.}
\label{fig:nmf_data}
\end{figure}

Here we assume a generic situation where signals from multiple sources are added as is common in non- or weakly interacting systems.
Examples in astronomy are mixture spectra from multiple stellar populations or spatial distributions of multiple galaxy types in galaxy clusters.

The most conservative option is the canonical non-negative matrix factorization (NMF), i.e. the parameterization and loss function from \autoref{eq:mf} with the penalty function
\begin{equation}
\iota_+(\x) = \begin{cases}0 & \mathrm{if}\ x_i \geq 0\ \forall i,\\\infty & \mathrm{else}\end{cases}
\end{equation}
for both matrix factors $\tA$ and $\tS$.
It provides a prototypical example of an efficient proximal operator $\prox_+(\x)=\max(0, \x)$, i.e. an element-wise thresholding operator.
The test data has $C=100$ observations with Gaussian i.i.d. noise of a mixture model of $K=3$ sinosoidal components $\in\mathbb{R}^{50}$ and is shown in \autoref{fig:nmf_data}.

We also run a variant of NMF, dubbed MixMF, that is additionally constrained to impose the mixture-model characteristic of these data, i.e. $\sum_k\tA_{ck} =1\ \forall c$.
The correspondent proximal operator is the projection operator onto the simplex, $\prox_{\mathrm{unity}}(\x) = |\x| / \sum |x_i|$, and is applied to every row $\tA_c$ for $c = \lbrace 1,\dots,C\rbrace$ to normalize the contributions of all components.

\begin{figure}
    \includegraphics[width=\linewidth]{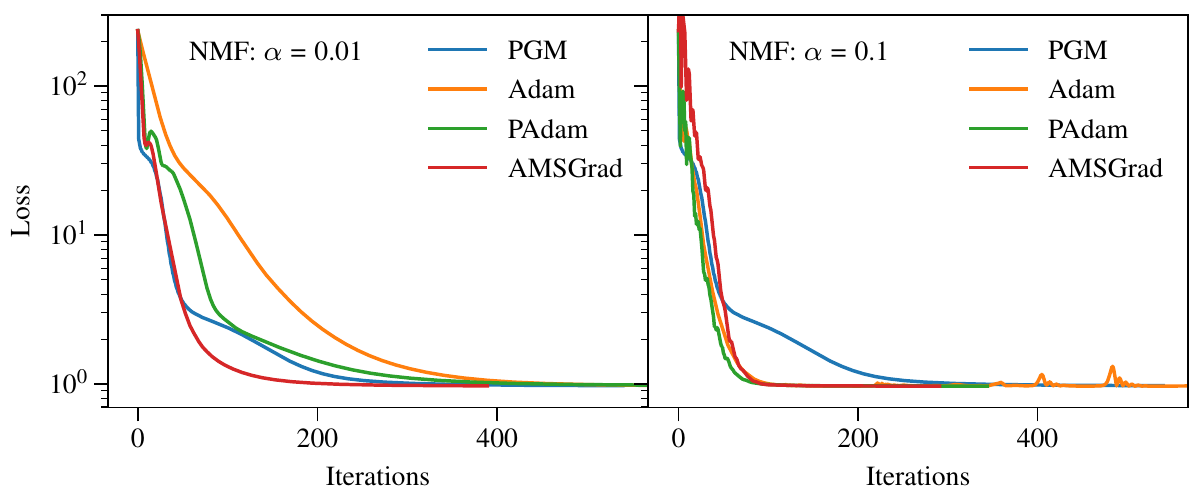}\\
    \includegraphics[width=\linewidth]{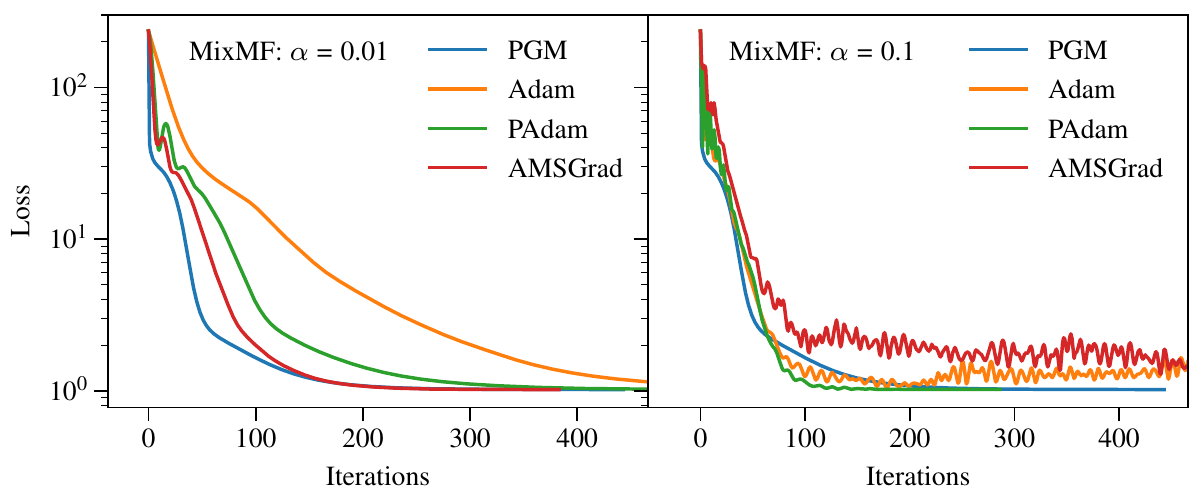}
    \caption{Loss for the NMF and MixMF problems of PGM and \adaprox\ with different step sizes and adaptive optimization schemes from \autoref{tab1}. Following the recommendation by \citet{Kingma2015-pq}, we set $\beta_1=0.9$, $\beta_2=0.999$, and $\epsilon=10^{-8}$; for \padam, we set $p=0.125$ as recommended by \citet{Chen2018-wu}. With fixed step sizes $\alpha$, \amsgrad\ and \adamx\ behave identically, the latter is thus not shown. Solutions are considered converged when the relative deviation of $\tA$ and $\tS$ between subsequent iterations is $<10^{-4}$.
}
    \label{fig:nmf_loss}
\end{figure}

We compare the performance in terms of final loss and number of evaluations and proximal evaluations for PGM and \adaprox with the adaptive schemes listed in \autoref{tab1}.
The initial values for $\tA$ and $\tS$ are drawn from $\mathcal{U}(0,1)$.
For PGM, we compute the analytic Lipschitz constants at every step.
For \adaprox, we choose the step sizes by considering the amplitude of the elements of $\tA$ and $\tS$, which are of order unity.
In the first run, we set them conservatively to $\alpha=0.01$, in the second run more aggressively to $\alpha=0.1$.
In both runs the step sizes are kept constant.
The results are shown in \autoref{fig:nmf_loss} and summarized in \autoref{tab2}.

\begin{table*}[tp]
\caption{Performance for the NMF (\emph{top}) and MixMF (\emph{bottom}) problem of PGM and \adaprox\ with different adaptive optimization schemes. See \autoref{fig:nmf_loss} for details. We list the number of iterations and the average number of proximal sub-iterations per iteration for $(\tA,\tS)$, respectively.}
\label{tab2}
{\small
\begin{tabular}{lcccccc}
\toprule
NMF & \multicolumn{3}{c}{$\alpha=0.01$} & \multicolumn{3}{c}{$\alpha=0.1$} \\
 \cmidrule(r){2-4} \cmidrule(r){5-7}
 Name & Final Loss & Iterations & Sub-Iterations & Final Loss & Iterations & Sub-Iterations\\
\midrule
PGM & 0.97261 & 541 & (1,1) & 0.97261 & 541 & (1,1)\\
\adaprox-\adam & 0.97121 & 663 & (1.89, 1.99) & \bf{0.96585} & 677 & (1.82, 1.98)\\
\adaprox-\padam & 0.97811 & 600 & (1.92, 1.92) & 0.96722 & 352 & (1.98, 1.98)\\
\adaprox-\amsgrad & \bf{0.96928} & \bf{405} & (1.97, 1.96) & 0.96645 & \bf{299} & (2.00, 2.00)\\
\midrule
MixMF & \multicolumn{3}{c}{$\alpha=0.01$} & \multicolumn{3}{c}{$\alpha=0.1$} \\
 \cmidrule(r){2-4} \cmidrule(r){5-7}Name & Final Loss & Iterations & Sub-Iterations & Final Loss & Iterations & Sub-Iterations\\
\midrule
PGM & 1.0193 & 444 & (1,1) & \bf{1.0193} & 444  & (1,1)\\
\adaprox-\adam & 1.0208 & 756 & (12.2, 1.99) & 1.3738 & 1000* & (16.5, 1.74)\\
\adaprox-\padam & 1.0208 & 528 & (2.74, 1.94) & 1.0227 & \bf{286} & (3.07, 1.87)\\
\adaprox-\amsgrad  & \bf{1.0191} & \bf{375} & (9.95, 1.99) & 1.3436 & 1000* & (11.1, 1.41)\\
\bottomrule
\end{tabular}
}
\\[0.2em]
{\small\it * indicates non-convergence after 1000 steps}
\end{table*}

Even with the conservative step sizes $\alpha=0.01$, \adaprox-\amsgrad\ outperforms PGM in terms of final loss and number of iterations for both problems.
For $\alpha=0.1$, every adaptive scheme outperforms PGM on the NMF problem, but they all show mild to prominent oscillations on the MixMF problem.
We find empirically that for \amsgrad\ and \padam, but not for \adam, this behavior can be mitigated by reducing the moving average parameters.
Values of $\beta_1\approx 0.5$ and $\beta_2\approx 0.8$ appear useful compromises between maintaining memory of previous gradients and adjusting to the newly constrained state.

Unsurprisingly, the solvers for the MixMF problem require more proximal sub-iterations for $\tA$ than for $\tS$ or for either in the NMF problem, where $\lesssim 2$ iterations of \autoref{eq:subiter} is sufficient.
That low number is due to the per-element thresholding of $\prox_+$.
If an element $x_i<0$, it will be projected to 0 on the first interation of $\prox_+$.
Because $\tH$ is diagonal, no other element is affected in the second iteration, the thesholding operation comes to the same result, and the sub-problem terminates after two iterations.
The mixture-model constraint, on the other hand, affects all elements of $\tA$ and therefore requires multiple passes to converge to the optimal solution.
It is intriguing that  \adaprox-\padam\ requires fewer calls of $\prox_{\mathrm{unity}+}$.
We do not have an explanation for this behavior.

Repeating the tests with different random seeds we establish these algorithm traits.
\begin{itemize}
\item Differences in the final loss between PGM and \adaprox\ are small.
Adaptive schemes can reach convergence in fewer iterations.
\item \adaprox-\adam\ shows good performance with small step sizes, but is the most unstable scheme overall, possibly related to the concerns raised by \citet{Reddi2018-me}.
\item \adaprox-\amsgrad\ shows some instability with larger step sizes; reducing $\beta_1$ and $\beta_2$ is beneficial.
\item \adaprox-\padam, using $p=[0.1, 0.25]$,  shows fast initial drops in the loss for large step sizes and converges quickly but to a slightly inferior final loss.
\end{itemize}

\subsection{Multi-band Source Separation}

We present a simplified test case that captures the main characteristic of source separation in astronomical imaging data observed in multiple optical filter bands (see \citet{Melchior2018-hb} for a full implementation).
The data set is comprised of $30\times30$ pixel images, observed in $C=5$ different filter bands, and affected by filter-dependent uncorrelated Gaussian background noise.
Using the definitions from \autoref{eq:fusion-mf}, and interpreting every filter as an independent observation $l$, we set $\Sigma_l^{-1} = \sigma_l^{-2} \mathbf{1}$.
The degradation operator $\mathsf{P}_l$ amounts to a simple projection of the hyperspectral data cube to a single observed filter band.
For the sake of simplicity and unlike actual observations, no convolution degrades the spatial resolution of the images.

\begin{figure}
    \includegraphics[width=\linewidth]{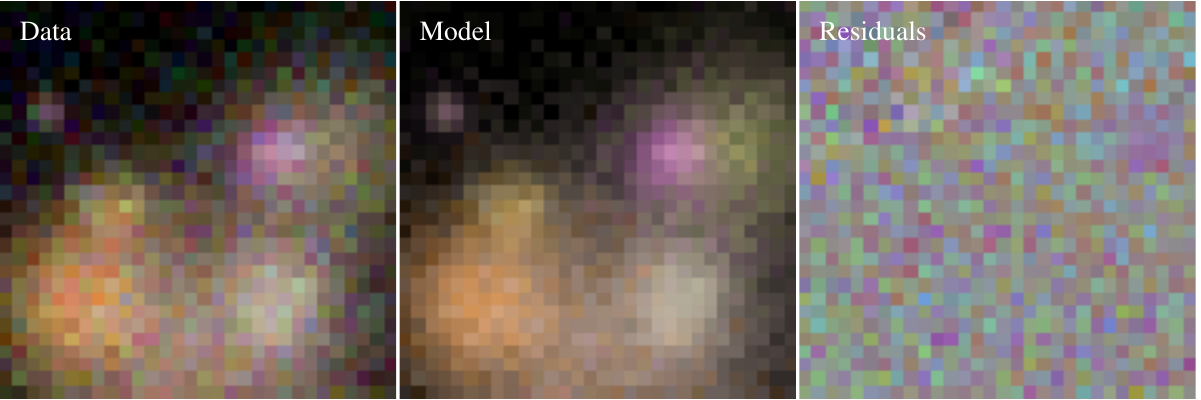}
    \caption{Astronomical source separation example. A false-color composite of a 5-band image data cube (\emph{left}) comprising $K=7$ circular Gaussian sources with band-dependent Gaussian additive noise; the CMF model of this scene from \adaprox-\padam\ (\emph{center}); individual components shown in \autoref{fig:astro_comp}) and its residuals (\emph{right}). The left and center panel use an inverse hyperbolic sine stretch \citep{Lupton2004-wv} to increase the dynamical range; the right panel uses a linear stretch. The test data are publicly available in the code repository.}
    \label{fig:astro_data}
\end{figure}

We distribute $K=7$ two-dimensional circular Gaussian-shaped sources randomly in the image, with sizes $\sigma_k$ ranging from 1 to 10 pixels.
Their integrated fluxes scale with the size, $F_k\propto \sigma_k^2$, as is approximately observed for astronomical sources.
The example multi-band image is shown in the left panel of \autoref{fig:astro_data}.

Since all astronomical sources are expected to be emitters of light, we impose a non-negativity constraint on $\tA$ and $\tS$.
We also add an $\ell_0$ penalty for $\tS$, whose proximal operator is the element-wise hard thresholding operator
\begin{equation}
\prox_{\lambda\, \ell_0}(\x) = \begin{cases} \x\ \ \mathrm{if}\ |\x| > \lambda\\0\ \ \mathrm{else}\end{cases}
\end{equation}
and then normalize the sum of the pixels with $\prox_{\mathrm{unity},+}$.
We initialize the individual components with circular Gaussians, whose centers and sizes are randomized by up to $\sigma_k/4$ and 50\%, respectively; their per-band amplitudes are taken from the noisy images at the assumed center.
This approach mimics a data processing pipeline that performs the initial object detection and characterization to warm-start the source separation method.

For PGM, we again compute the analytic Lipschitz constants at every iteration to set the step sizes.
For \adaprox, we follow the general logic of \autoref{sec:nmf} and set them relative to their typical amplitude.
As the spatial distributions are normalized to unity, their mean amplitude is $\approx 10^{-3}$, and we decide on a more conservative setting by adjusting the step sizes to 1\%, i.e. $\alpha^{(\tS)}=10^{-5}$.
The per-band amplitudes $\tA_k$, however, are different by a factor of $\approx100$ between the brightest and the faintest source.
Unlike in PGM, we are free to set them differently for every component and chose $\alpha^{(\tA_k)}=\tfrac{0.1}{C}\sum_c \tA_{ck}$, reflecting our expectation that the initial amplitudes can have errors on the order of 10\%.

\begin{figure}
    \includegraphics[width=\linewidth]{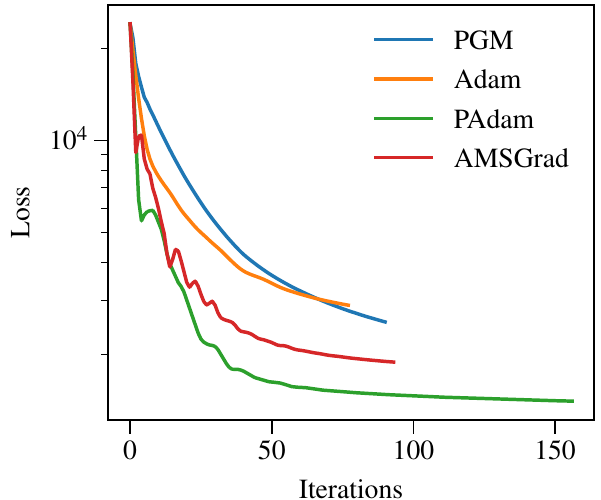}
    \caption{Loss for the astronomy CMF problem of PGM and \adaprox\ with different adaptive optimization schemes from \autoref{tab1}.
    Following the recommendation by \citet{Kingma2015-pq}, we set $\beta_1=0.9$, $\beta_2=0.999$, and $\epsilon=10^{-8}$; for \padam, we set $p=0.45$, the best-performing value for this problem (see details in the text).
    Solutions are considered converged when the relative deviation of $\tA$ and $\tS$ between subsequent iterations is $<10^{-3}$.}
    \label{fig:astro_loss}
\end{figure}

The resulting losses are shown in \autoref{fig:astro_loss} and summarized in \autoref{tab3}.
It is evident that \adaprox\ can match or outperform PGM in terms of the final loss within a similar number of iterations.
\adaprox-\padam\ yields the best result, albeit with a slower convergence, but only after adjusting $p=0.45$.
With the recommended $p\approx0.125$, the first few steps move far away from the initial $\tA$ and $\tS$, which means the solver largely ignores the reasonable starting positions.
At $p=0.5$, \padam\ is identical to \amsgrad.
Intermediate values of $p$ appear to compromise between rapid initial improvement of the loss with \padam\ and the robust performance that characterizes \amsgrad\ at smaller step sizes, in accordance with our observations in \autoref{sec:nmf}.
We note, however, that the step sizes are only given in units of the parameter if $p=0.5$ so that the gradient amplitude is cancelled by the term $\phi/\psi$.

With the given constraints, \adaprox\ requires only 1 to 2 proximal sub-iterations.
By avoiding the computation of the spectral norm for the Lipschitz constants in PGM, \adaprox\ exhibits much lower runtimes, which more than outweighs the computational cost of the adaptive schemes and the extra evaluations of the proximal operators.

\begin{table}[tp]
\setlength\tabcolsep{1.7pt}
\caption{Performance for the astronomy CMF problem problem of PGM and \adaprox. See \autoref{fig:astro_loss} for details. We list the number of iterations and the number of proximal sub-iterations per iteration for $(\tA,\tS)$, respectively. The runtime is for a single CPU on a recent Apple MacBook Pro.}
\label{tab3}
{\small
\begin{tabular}{lcccc}
\toprule
Name & Loss & Iterations & Sub-Iterations & Runtime [s]\\
\midrule
PGM & 2538.4 & 91 & (1,1) & 6.3 \\
\adaprox-\adam & 2884.8 & \bf{78} & (1.97, 2.00) & 0.17\\
\adaprox-\padam & \bf{1398.2} & 167 & (1.19, 2.13) & 0.38\\
\adaprox-\amsgrad & 1883.9 & 94 & (1.51, 2.00) & \bf{0.12}\\
\bottomrule
\end{tabular}
}
\end{table}

\begin{figure}
    \includegraphics[width=\linewidth]{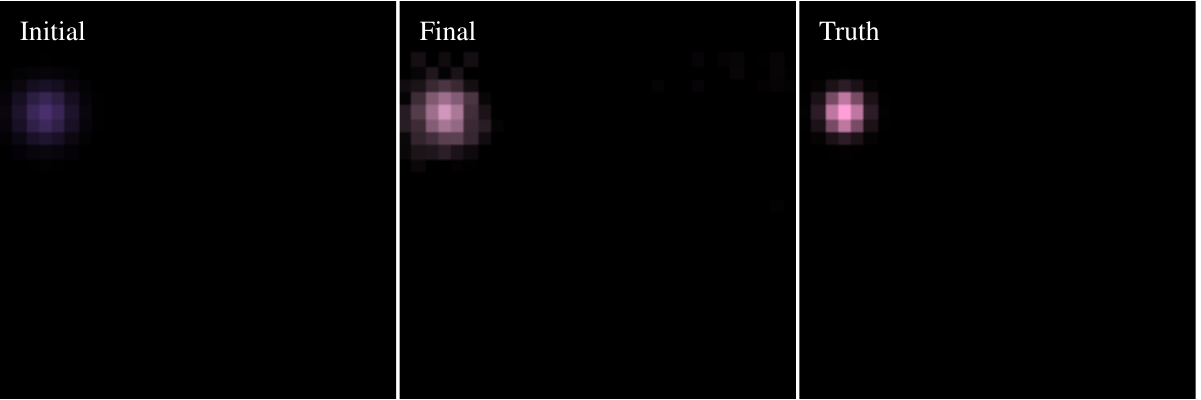}\\
    \includegraphics[width=\linewidth]{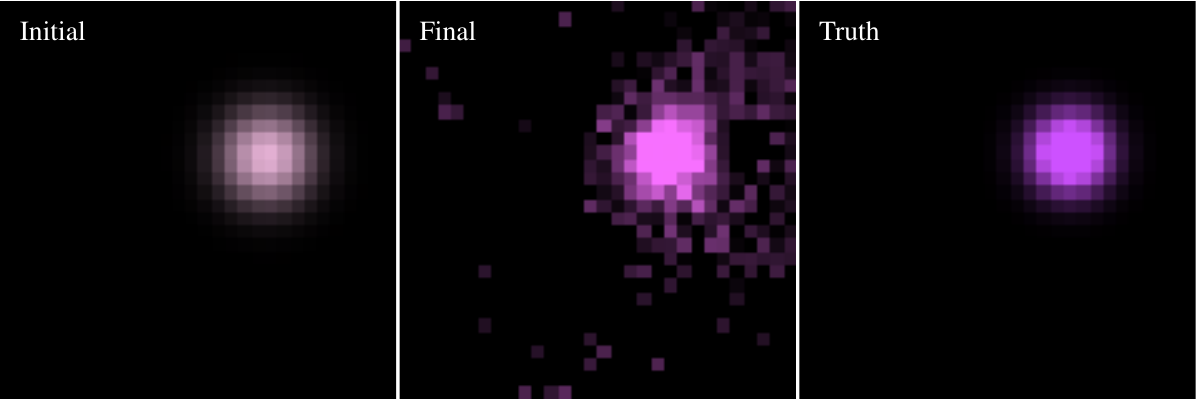}\\
    \includegraphics[width=\linewidth]{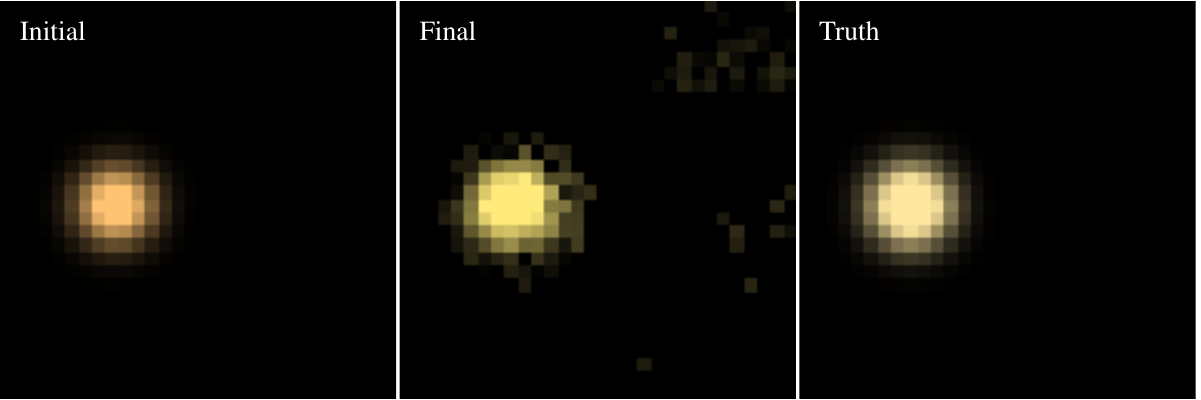}\\
    \includegraphics[width=\linewidth]{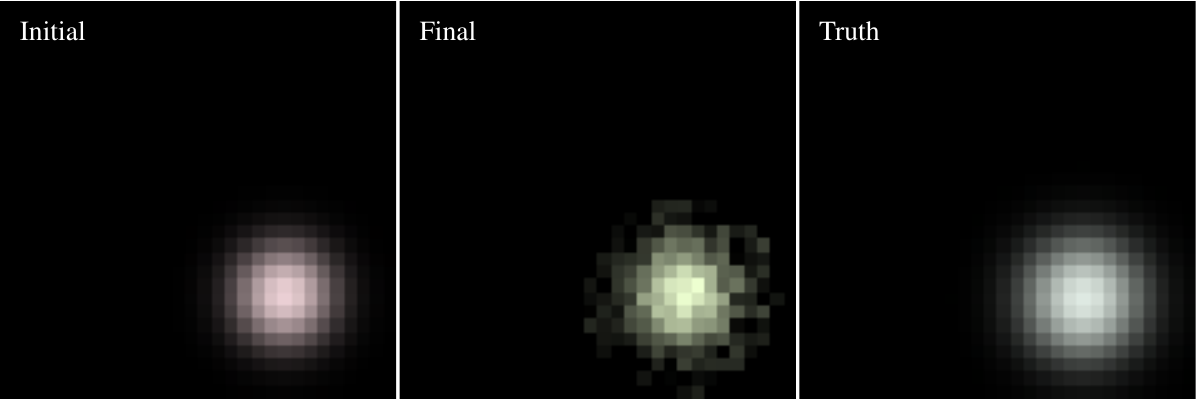}\\
    \includegraphics[width=\linewidth]{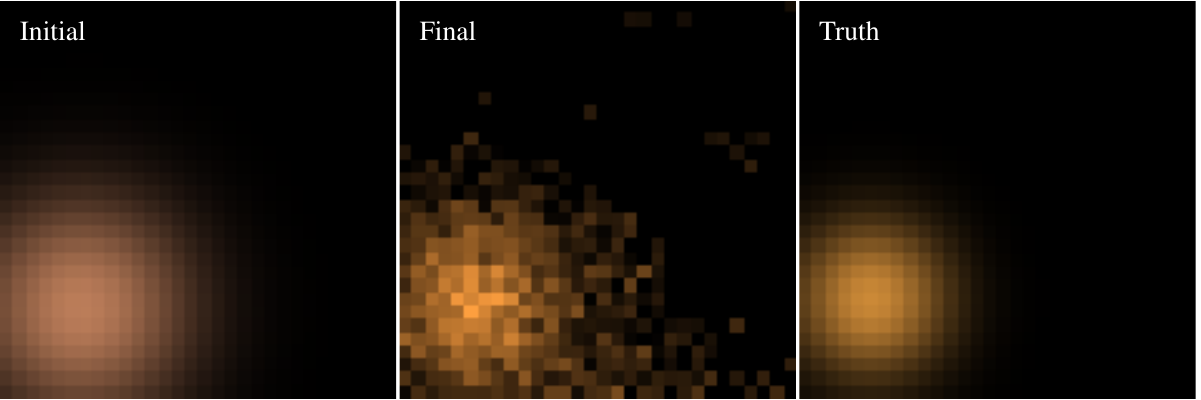}\\
    \includegraphics[width=\linewidth]{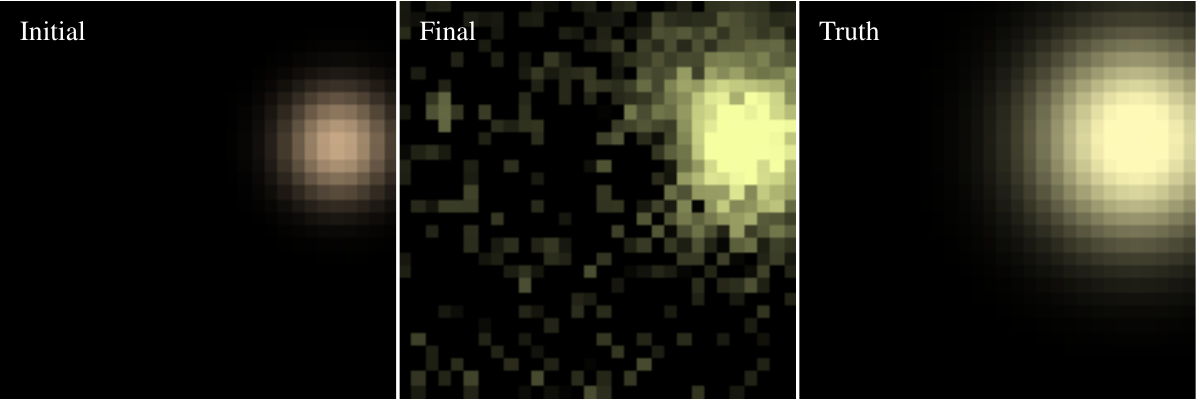}\\
    \includegraphics[width=\linewidth]{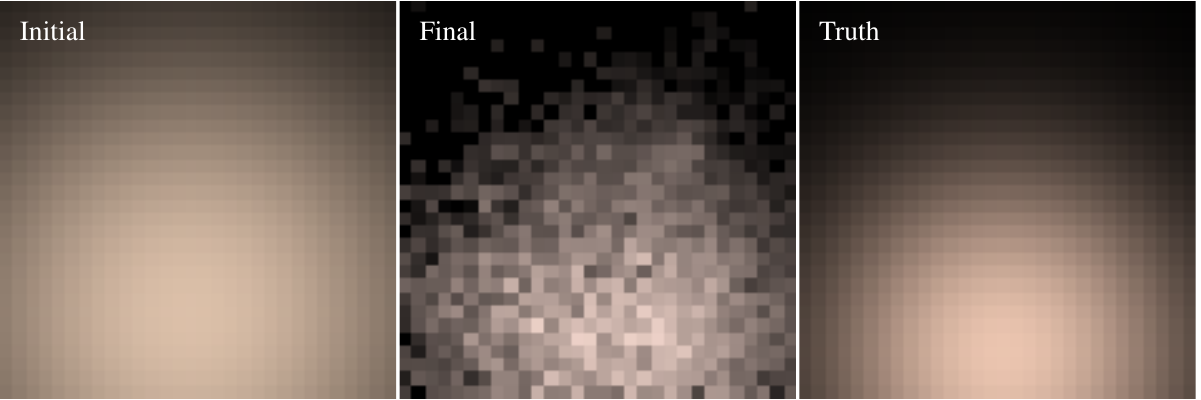}
    \caption{Individual components of the model shown in the middle panel of \autoref{fig:astro_data}. Each component is initialized with the best-fitting Gaussian to the subregion in the image. The images use an inverse hyperbolic sine stretch adjusted for each component.}
    \label{fig:astro_comp}
\end{figure}

The visual inspection of the individual components (\autoref{fig:astro_comp}) of the best-fitting \adaprox-\padam\ model confirms that the colors and shapes improve noticeably from well-chosen initial parameters.
The spatial distributions reveal the impact of noise but also the effect of the $\ell_0$ penalty, which promotes configurations with few non-zero pixels.

\section{Summary}
\label{sec:conclusions}

We present an adaptive proximal gradient method, \adaprox, which enables constrained convex optimization using the gradient updates of the recently proposed unconstrained me\-thod \adam\ and its variants \amsgrad, \adamx\  and \padam.
We solve the arising variable-metric proximal iteration by ordinary proximal gradient sub-iterations.
The scheme is applicable to arbitrary proxable penalty functions.
The cost of our proposed method arises from the need to compute and store moving averages of the first and second moment of the gradient of $f$ as well as multiple computations of the proximal mapping.
Its benefits stem from adjusting the effective learning rates for every parameters and from avoiding the computation of the Lipschitz constant, traditionally required for the proximal gradient method.
\adaprox\ is thus beneficial in cases when the Lipschitz constants cannot be calculated analytically or efficiently, e.g. for non-linear models or in signal processing problems with complicated observation designs, and when $f$ is too expensive to evaluate to determine $L$ through line searches.

We demonstrate in three variants of constrained matrix factorization problem that \adaprox, in particular with the \amsgrad\ and \padam\ schemes, outperforms PGM in terms of final loss, number of iterations, and runtime.
\adaprox\ requires that step sizes for each parameter are set in advance.
We find that relative step sizes on the order of 1\% to 10\% of the typical amplitude of the parameters work well in practice.

The {\sf python} implementation of the algorithms presented here are available as an open-source package at \url{https://github.com/pmelchior/proxmin}.

\section*{Acknowledgements}
We gratefully acknowledge partial support from NASA grant NNX15AJ78G.

For \software{proxmin} and the preparation of this manuscript we made use of the following software packages: \software{SciPy} \citep{2020SciPy-NMeth}, \software{numpy} \citep{numpy}, and \software{matplotlib} \citep{matplotlib}.

\bibliographystyle{spbasic}
\bibliography{references.bib}

\end{document}